\documentclass[oneside,10pt]
{amsart}          
\usepackage{graphicx}
\usepackage{amsfonts,
amsmath,latexsym,amssymb} 
\usepackage{amsthm}                
\usepackage{calrsfs}
\usepackage{enumerate}
\usepackage{enumitem}
\usepackage{url}
\usepackage[final]
{showkeys}

\usepackage{hyperref}

\theoremstyle{definition}

\newlist{thmRoman}{enumerate}{1}
\setlist[thmRoman]{%
  label=\normalfont
  \small(\Roman*), 
  ref=(\Roman*),                                    
  font=\normalfont
  \small,           
  align=left, leftmargin=*,                         
}

\newcommand{\al}{\alpha}

\newcommand{\ep}{\varepsilon}

\renewcommand{\Psi}{\overline{\Phi}}
\renewcommand{\P}{\operatorname{\mathsf{P}}} 

\newcommand{\F}
{\Sigma}

\renewcommand{\ge}{\geqslant}
\begin{document}

\title{Comment on ``Worst-case Nonparametric Bounds for the Student T-statistic'', arXiv:2508.13226
}


\author{Iosif Pinelis}

\address{Department of Mathematical Sciences\\
Michigan Technological University\\
Houghton, Michigan 49931, USA\\
}

%



\date{5 August 2015}                               

        



\begin{abstract}
Concerning Version~1 of ``Worst-case Nonparametric Bounds for the Student T-statistic'', arXiv:2508.13226: The main result there is incorrect. Concerning Version~2 of arXiv:2508.13226: At least the proof of the main result there is incorrect.  
\end{abstract}

\maketitle







The main result of the manuscript in question is Theorem~1, which claims that for each natural $n$ and each real $t>0$ 
\begin{equation*}
	\sup_{w\in W_n} \P(S_n(w)\ge t)=\max_{k\in[n]}\P(S_n(w^{(k)})\ge t),
\end{equation*}
where $W_n$ is the set of all vectors $w=(w_1,\dots,w_n)\in[0,\infty)^n$ such that $\sum_{j\in[n]}w_j^2=1$, $S_n(w):=\sum_{j\in[n]}w_j\ep_j$, $[n]:=\{1,\dots,n\}$, the $\ep_j$'s are independent Rademacher random variables (r.v.'s), and $w^{(k)}:=\frac1{\sqrt k}\,(1,\dots,1,0,\dots,0)\in W_n$. 

This claim contradicts \cite[Proposition~1]{sup-tails-radem-published}. In fact, this claim contradicts even the narrower result of \cite{zhubr_published}. There are no references in the manuscript in question to these previous results, even though it is clear from the correspondence that the author of the manuscript was aware of these results in \cite{sup-tails-radem-published,zhubr_published}. 

``The key observation'' in the first proof of Lemma 1 (which is described in the manuscript as ``The key insight underlying Theorem 1'') is about comparing so-called peakedness. For this, the author refers to a paper by Proschan -- which, however, deals with absolutely continuous distributions, whereas the distribution of $S_n(w)$ is of course discrete. 

The second proof of Lemma 1 seems to involve conditioning on $\ep_3,\dots,\ep_k$ (with $k$ somehow replacing $n$), but the notation used in that proof is for unconditional probabilities. Worst of all, inequality (12) in the manuscript fails to hold when, for instance, $\al=0$, $u=\frac38$, $v=\frac58$, and $t=\frac18$. 

\medskip 
\hrule 
\medskip 

After the previous version of this comment was posted, the author of manuscript arXiv:2508.13226 has alerted me to Version~2 of the manuscript. In this version, the author added a reference to \cite{sup-tails-radem-published} but not to \cite{zhubr_published}. 

The main result in the manuscript and the corresponding proof have now been changed significantly. Yet, even a quick glance at the revised manuscript still shows issues of different degrees of severity. 

One issue is that the maximum in (2) does not have to be attained in general. 

The proof of the main result (that is, of Theorem~1) is less than three lines -- basically, just a reference to Lemma~1. I do not understand this proof, even assuming that Lemma~1 is correct. 

As for Lemma~1: 
Taken literally, the first sentence after formula~(6) in the proof of the lemma does not make sense, since neither side of (6) depends on $\theta$. Also, if $\ep_j=1$ then the right-hand side of (6) will always be, not positive, but negative -- given the condition $w_i>w_j>0$.

There are a number of other issues. 

Perhaps most importantly, the proof of Lemma~1 is based on the claim that the c.d.f.\ of a symmetric discrete, generally nonzero r.v.\ $T$ is concave on $[0,\infty)$. However, such a c.d.f.\ is not continuous on $(0,\infty)$ and therefore cannot possibly be concave on $[0,\infty)$. 

The false conclusion on the concavity is due to the mentioned author's misunderstanding of the notion of unimodality and related facts considered in Ref.~[15] in the manuscript in question and used in the proof of Lemma~1 there. First here, the notion of unimodality considered in Ref.\ [15] pertains only to  distributions that are absolutely continuous to the left and to the right of some point, whereas the distributions considered in the proof of mentioned Lemma~1 are discrete and generally non-degenerate. Second, it is not true in general that the convolution of unimodal distributions is unimodal (see e.g.\ \cite[Appendix II]{gned-kolm}) -- which is why the notion of \emph{strong} unimodality is indtroduced (and characterized) in mentioned Ref.~[15]. 


\bibliographystyle{abbrv}


\bibliography{C:/Users/ipinelis/Documents/pCloudSync/mtu_pCloud_02-02-17/bib_files/citations04-02-21}


\end{document}